\documentclass[12pt]{article}
\usepackage{amsmath,amssymb,amscd,amsthm,latexsym}

\font\block=msbm10
\def\C{\hbox{\block\char'0103}}
\def\Z{\hbox{\block\char'0132}}
\def\+{\hbox{\block\char'0156}}

\font\gotic=eufm10
\def\g{\hbox{\gotic\char'0147}}
\def\l{\hbox{\gotic\char'0154}}
\def\s{\hbox{\gotic\char'0163}}
\def\p{\hbox{\gotic\char'0160}}
\def\o{\hbox{\gotic\char'0166}}

\def\u{\hbox{\gotic\char'0165}}
\def\c{\hbox{\gotic\char'0143}}

\font\new=eusm10
\def\A{\hbox{\new\char'0101}}
\def\P{\hbox{\new\char'0120}}

\def\W{\hbox{\new\char'0127}}

\font\yes=msam10
\def\Q{\hbox{\yes\char'0003}}

\begin{document}

\noindent
\centerline{\Large On matrix realizations of the Lie superalgebra $D(2, 1 ; \alpha)$}

\vskip 0.5in
{\large Elena Poletaeva}
\vskip 0.1in

{\it Department of Mathematics,}

{\it University of Texas-Pan American,}

{\it Edinburg, TX 78539}

{\it Electronic mail:} elenap$@$utpa.edu

\vskip 0.5in
\noindent
{\footnotesize \noindent {\bf Abstract.}
We obtain a realization of the Lie
superalgebra $D(2, 1 ; \alpha)$ in differential operators on the supercircle $S^{1|2}$ and
in $4\times 4$ matrices over a Weyl algebra.
A contraction of $D(2, 1 ; \alpha)$
is isomorphic to the universal central extension $\hat{\p\s\l}(2|2)$ of $\p\s\l(2|2)$.
We realize it in $4\times 4$ matrices over the associative algebra of pseudodifferential operators on $S^1$.
Correspondingly, there exists a three--parameter family of irreducible representations of
$\hat{\p\s\l}(2|2)$ in a $(2|2)$--dimensional complex superspace.}

\vskip 0.1in
\noindent
{\it JGP SC:} Lie superalgebras
\vskip 0.1in
\noindent

\vskip 0.5in
\noindent
{\bf 1. Introduction}

\vskip 0.3in

Recall that $D(2, 1 ; \alpha)$, where  $\alpha \in \C\backslash
\lbrace 0, -1 \rbrace$, is a one-parameter family of classical simple exceptional Lie superalgebras of dimension 17 [1].
The bosonic part of $D(2, 1 ; \alpha)$ is
$\s\l(2)\oplus \s\l(2)\oplus \s\l(2)$, and the action of
$D(2, 1 ; \alpha)_{\bar 0}$ on $D(2, 1 ; \alpha)_{\bar 1}$
is the  product of two-dimensional representations.
This family has a very interesting relation to the simple classical Lie superalgebra $\p\s\l(2|2)$,
which is the only basic classical Lie superalgebra having a nontrivial universal central extension  with three central elements [2]. Note that this central extension can be obtained as a contraction of $D(2, 1 ; \alpha)$ when
$\alpha\rightarrow -1$. The Lie superalgebra $\p\s\l(2|2)$ and its central extensions play an important r{\^o}le in the AdS/CFT correspondence [3--6].

In [7] M. Scheunert gave an equivalent description of the superalgebra $D(2, 1 ; \alpha)$ as the family of superalgebras
$\Gamma (\sigma_1, \sigma_2,  \sigma_3)$, where $\sigma_i$ are nonzero complex numbers
such that $\sigma_1+\sigma_2+\sigma_3 = 0$.
Note that $\Gamma (\sigma_1, \sigma_2, \sigma_3)\cong D(2, 1 ; \alpha)$, where
$\alpha = \sigma_1/{\sigma_2}$. We use the notation of [7].

In this work we consider the standard embedding of
$\Gamma (2, -1 - \alpha,  \alpha - 1)$,
where $\alpha \in \C$,
into the Poisson superalgebra of {\it differential} operators on the
supercircle $S^{1|2}$ with two odd variables and into its deformation.
This allows us to realize $\Gamma (2, -1 - \alpha,  \alpha - 1)$ in $4\times 4$ matrices over the Weyl algebra
$\W = \sum_{i\geq 0} \C[t, t^{-1}]d^i$, where
$d = {\partial\over {\partial t}}$.
This realization differs from the realization that we obtained in [8, 9], where we essentially used
{\it pseudodifferential} operators on $S^{1|2}$.

Note that in [10, 11] we realized the superconformal algebras  $\hat{K}'(4)$ and $CK_6$ in matrices over $\W$
of size 4 and 8, respectively.
A superconformal algebra naturally has a $\Z$-grading (see [12]).
The zeroth part of the $\Z$-grading of $\hat{K}'(4)$ is isomorphic to
a central extension of $\p\s\l(2|2)$ with two central elements.
The zeroth part of the $\Z$-grading of $CK_6$ is isomorphic to
the universal central extension $\hat{P}(3)$ of the classical simple Lie superalgebra $P(3)$ (see [1] for notation).
C. Martinez and E. I. Zelmanov obtained the analogous realization of $\hat{P}(3)$ in matrices of size 8 over $\W$
by a different approach using Jordan superalgebras; see [13].

In [14], C. Martinez, I. Shestakov and E. Zelmanov
considered Cheng-Kac Jordan superalgebras $\hbox{CK}(Z, D)$, where $Z$ is an associative commutative algebra with a derivation $D: Z\longrightarrow Z$. They obtained an embedding of
$\hbox{CK}(Z, D)$ into the  $4\times 4$ matrices over the algebra $\W$ of differential operators on $Z$ with the $\Z/2\Z$-gradation. This embedding extends the King-McCrimmon embedding of the Kantor double of the vector type bracket $\lbrace a, b\rbrace = D(a)b - aD(b)$ into
$2\times 2$ matrices over $\W$; see [15, 16].
The Cheng-Kac Lie superalgebra  introduced in [17] and denoted by  $CK_6$ is isomorphic to the Tits-Kantor-Koecher construction of the Jordan superalgebra
$\hbox{CK}(\C[t, t^{-1}], d/dt)$. In [14], the authors obtained, in particular,
 an embedding of the
Cheng-Kac Lie superalgebra into a superalgebra of $8\times 8$ matrices over the Weyl algebra $\W$.

The methods of [14] can be used to obtain an embedding of $\Gamma (\sigma_1, \sigma_2, \sigma_3)$.
On the Jordan level, the problem is equivalent to an embedding of four-dimensional Jordan superalgebras $D(\alpha)$ into $2\times 2$-matrices over $\W$. Let $R = \C[t]$ and $d = d/dt$. Consider the algebra
$\tilde{R} = R + Rw$, $w^2 = 1$  with the derivation $d$, and consider the vector type bracket
$\lbrace a, b\rbrace = d(a)b - ad(b)$.
Then the Kantor double $K(\tilde{R}, d) = \tilde{R} + \tilde{R}v$ contains
$\C +\C w + \C v +\C[(1-\alpha)tw + (1 + \alpha)t]v\cong D(\alpha)$.

K. McCrimmon proved in [18] that a Kantor double of a vector type bracket can be embedded into
$2\times 2$ matrices over the corresponding Weyl algebra.
Then there exists an embedding of $D(\alpha)$ into $2\times 2$-matrices over $\W \oplus \W$, and also an embedding of $D(\alpha)$ into $2\times 2$ matrices over $\W$, since $D(\alpha)$ is simple.

A contraction of $\Gamma (2, -1 - \alpha,  \alpha - 1)$ when $\alpha \rightarrow \pm 1$ is isomorphic to the universal central extension $\hat{\p\s\l}(2|2)$ of $\p\s\l(2|2)$.
We realize it in $4\times 4$ matrices over the associate algebra of pseudodifferential operators on $S^{1}$, which contains $\W$ as a subalgebra.  Correspondingly, we obtain a three-parameter family of irreducible representations of $\hat{\p\s\l}(2|2)$ in $(2|2)$--dimensional complex superspace.

\vskip 0.5in
\noindent
{\bf 2. Embedding of $\Gamma(\sigma_1, \sigma_2, \sigma_3)$ into the Poisson superalgebra on $S^{1|2}$}

\vskip 0.5in

Recall the definition of $\Gamma(\sigma_1, \sigma_2, \sigma_3)$ [7].
Let $\g = \g_{\bar{0}} \oplus \g_{\bar{1}}$ be a Lie superalgebra, where
$\g_{\bar{0}} = sp(\psi_1)\oplus sp(\psi_2)\oplus sp(\psi_3)$ and
$\g_{\bar{1}} = V_1\otimes V_2\otimes V_3$, where
$V_i$ are two-dimensional vector spaces, and
$\psi_i$ is a non-degenerate skew-symmetric form on $V_i$, $i = 1, 2, 3$.
A representation of $\g_{\bar{0}}$ on $\g_{\bar{1}}$ is the tensor product
of the standard representations of $sp(\psi_i)$ in $V_i$.
Consider the $sp(\psi_i)$ - invariant bilinear mapping
$$\P_i: V_i\times V_i \rightarrow sp(\psi_i), \quad i = 1, 2, 3,\eqno (1)$$
given by
$$\P_i(x_i, y_i)z_i = \psi_i(y_i, z_i)x_i - \psi_i(z_i, x_i)y_i\eqno (2)$$
for all $x_i, y_i, z_i\in V_i$.
Let $\P$ be a mapping
$$\P:\g_{\bar{1}}\times \g_{\bar{1}}\rightarrow \g_{\bar{0}}\eqno (3)$$
given by
\begin{equation*}
\begin{aligned}
&\P(x_1\otimes x_2\otimes x_3, y_1\otimes y_2\otimes y_3) =\\
&\sigma_1\psi_2(x_2, y_2)\psi_3(x_3, y_3)\P_1(x_1, y_1) + \\
&\sigma_2\psi_1(x_1, y_1)\psi_3(x_3, y_3)\P_2(x_2, y_2) + \\
&\sigma_3\psi_1(x_1, y_1)\psi_2(x_2, y_2)\P_3(x_3, y_3)
\end{aligned}
\tag{4}
\end{equation*}
for all $x_i, y_i \in V_i, i = 1, 2, 3$,
where $\sigma_1, \sigma_2, \sigma_3$ are some complex numbers.
The super-Jacobi identity is satisfied if and only if
$\sigma_1 + \sigma_2 + \sigma_3 = 0$. In this case $\g$ is denoted by
$\Gamma(\sigma_1, \sigma_2, \sigma_3)$.
Superalgebras $\Gamma(\sigma_1, \sigma_2, \sigma_3)$ and
$\Gamma(\sigma_1', \sigma_2', \sigma_3')$ are isomorphic if and only if there exists a nonzero element $k\in \C$ and a permutation $\pi$ of the set $\lbrace 1, 2, 3\rbrace$ such that
$$ \sigma_i' = k\cdot\sigma_{\pi i} \hbox{ for } i = 1, 2, 3.\eqno (5)$$
Superalgebras $\Gamma(\sigma_1, \sigma_2, \sigma_3)$ are simple if
and only if $\sigma_1, \sigma_2, \sigma_3$ are all different from zero.
Note that $\Gamma(\sigma_1, \sigma_2, \sigma_3)\cong D(2, 1 ; \alpha)$ (see [1]) where
$\alpha = \sigma_1/{\sigma_2}$.

 The {\it Poisson algebra $P$ of pseudodifferential operators on
the circle} is formed by the formal series
$$A(t, \tau) = \sum_{-\infty}^na_i(t)\tau^i,\eqno (6)$$
where $a_i(t)\in \C[t, t^{-1}]$,
and the even variable $\tau$ corresponds to $\partial_t$, see [19].
The Poisson bracket is defined as follows:
$$\lbrace A(t, \tau), B(t, \tau)\rbrace = \partial_{\tau}A(t, \tau)\partial_tB(t, \tau) - \partial_tA(t, \tau)\partial_{\tau}B(t, \tau).\eqno (7)$$
An associative algebra
$P_{\hbox{h}}$, where $\hbox{h}\in  (0, 1]$, is a deformation of $P$, see [20].
The multiplication in $P_{\hbox{h}}$ is given as follows:
$$A(t, \tau)\circ_{\hbox {h}}B(t, \tau) = \sum_{n\geq 0 }{{\hbox{h}}^n\over {n!}}
\partial_{\tau}^nA(t, \tau)\partial_t^nB(t, \tau).\eqno (8)$$
The Lie algebra structure on the vector space $P_{\hbox{h}}$ is given by
$$[A, B]_{\hbox{h}} = {1\over \hbox{h}}(A\circ_{\hbox{h}}B - B\circ_{\hbox{h}}A),\eqno (9)$$
and so
$$\hbox{lim}_{\hbox{h}\rightarrow 0}[A, B]_{\hbox{h}} = \lbrace A, B\rbrace.\eqno (10)$$
Let $\Lambda(2N)$ be the Grassmann algebra in $2N$ variables
$\xi_1, \ldots, \xi_N, \eta_1, \ldots, \eta_N $ with the parity $p(\xi_i) = p(\eta_i) = \bar{1}$.
The {\it Poisson superalgebra of pseudodifferential operators on $S^{1|N}$} is
$P(2N) = P\otimes \Lambda(2N)$. The Poisson bracket is defined as follows:
$$\lbrace A, B\rbrace = \partial_{\tau}A\partial_tB - \partial_tA\partial_{\tau}B +
(-1)^{p(A)+1}\sum_{i = 1}^N(\partial_{\xi_i}A\partial_{\eta_i}B + \partial_{\eta_i}A\partial_{\xi_i}B).\eqno (11)$$
Let $\Lambda_{\hbox{h}}(2N)$ be an associative superalgebra with generators
$\xi_1, \ldots, \xi_N, \eta_1, \ldots, \eta_N$ and relations
$$\xi_i\xi_j = - \xi_j\xi_i,\quad \eta_i\eta_j = -\eta_j\eta_i, \quad \eta_i\xi_j = \hbox{h}\delta_{i, j} - \xi_j\eta_i.\eqno (12)$$
Let $P_{\hbox{h}}(2N) = P_{\hbox{h}}\otimes\Lambda_{\hbox{h}}(2N)$ be a superalgebra with the product given by
$$(A_1\otimes X)(B_1\otimes Y) =
(A_1\circ_{\hbox{h}}B_1)\otimes (XY),\eqno (13)$$
 where
$A_1, B_1 \in P_{\hbox{h}}$ and $X, Y \in \Lambda_{\hbox{h}}(2N)$.
The Lie bracket of $A = A_1\otimes X$ and $B = B_1\otimes Y$  is
$$[A, B]_{\hbox{h}} = {1\over \hbox{h}}(AB - (-1)^{p(X)p(Y)}BA),\eqno (14)$$
and (10) holds.
$P_{\hbox{h}}(2N)$ is called the {\it Lie superalgebra of pseudodifferential operators on $S^{1|N}$};
see [10, 11].
We also consider the  subalgebras of {\it differential operators} $P^+\subset P$ and $P_{\hbox{h}}^+\subset P_{\hbox{h}}$ formed by
$\sum_{i\geq 0}^na_i(t)\tau^i,$ and their superanalogues $P^+(2N)$ and $P_{\hbox{h}}^+(2N)$.

\vskip 0.1in
\noindent
{\bf Proposition 2.1.} For each $\alpha\in\C$
there exists an embedding
$$\rho_{\alpha}: \Gamma (2, -1 - \alpha, \alpha - 1)\subset P^+(4). \eqno (15)$$
$\Gamma_{\alpha} = \rho_{\alpha}(\Gamma (2, -1 - \alpha, \alpha - 1))$ is
spanned by the following elements:
\begin{equation*}
\begin{aligned}
&E_{\alpha}^1 =  t^2, \quad
F_{\alpha}^1 =  \tau^2 - 2\alpha t^{-2}\xi_1\xi_2\eta_1\eta_2, \quad
H_{\alpha}^1 =  t\tau,\\
&E_{\alpha}^2 = \xi_1\xi_2,\quad
F_{\alpha}^2 = \eta_1\eta_2, \quad
H_{\alpha}^2 =  \xi_1\eta_1 + \xi_2\eta_2, \\
&E_{\alpha}^3 =  \xi_1\eta_2, \quad
F_{\alpha}^3 = \xi_2\eta_1,\quad
H_{\alpha}^3 =  \xi_1\eta_1 - \xi_2\eta_2,\\
&T_{\alpha}^1 =  t\eta_1, \quad
T_{\alpha}^2 =  t\eta_2 ,\quad
T_{\alpha}^3 =  t\xi_1, \quad
T_{\alpha}^4 =  t\xi_2,\\
&D_{\alpha}^1 =  \tau\xi_1 + \alpha t^{-1}\xi_1\xi_2\eta_2,\quad
D_{\alpha}^2 = \tau\xi_2 - \alpha t^{-1}\xi_1\xi_2\eta_1,\\
&D_{\alpha}^3 =  \tau\eta_1 + \alpha t^{-1}\xi_2\eta_1\eta_2,\quad
D_{\alpha}^4 =  \tau\eta_2 - \alpha t^{-1}\xi_1\eta_1\eta_2.
\end{aligned}
\tag{16}
\end{equation*}

\noindent
{\it Proof.}
Note that if $\alpha = 0$, then
$\Gamma (2, -1, - 1)\cong \s\p\o (2|4)$, and $\rho_{\alpha}$ is the standard
embedding of $\s\p\o (2|4)$ into $P^+(4)$.

Let
\begin{equation*}
\begin{aligned}
&V_1 =\hbox{Span}(e_1, e_2) , \quad
V_2   = \hbox{Span}(f_1, f_2), \quad
V_3 = \hbox{Span}(h_1, h_2),
\end{aligned}
\tag{17}
\end{equation*}
and
\begin{equation*}
\begin{aligned}
&\psi_1(e_1, e_2) = - \psi_1 (e_2, e_1) = 1,\\
& \psi_2(f_1, f_2) = - \psi_2(f_2, f_1) = 1,\\
&\psi_3(h_1, h_2) = - \psi_3(h_2, h_1) = 1.
\end{aligned}
\tag{18}
\end{equation*}
Explicitly an embedding $\rho_{\alpha}$ is given as follows:
\begin{equation*}
\begin{aligned}
&\rho_{\alpha}(\P_1(e_1, e_1)) = -E_{\alpha}^1, \quad
\rho_{\alpha}(\P_1(e_2, e_2)) = -F_{\alpha}^1,\quad
\rho_{\alpha}(\P_1(e_1, e_2)) = -H_{\alpha}^1,\\
&\rho_{\alpha}(\P_2(f_1, f_1)) = - 2F_{\alpha}^2, \quad
\rho_{\alpha}(\P_2(f_2, f_2)) = - 2E_{\alpha}^2,\quad
\rho_{\alpha}(\P_2(f_1, f_2)) =  H_{\alpha}^2,\\
&\rho_{\alpha}(\P_3(h_1, h_1)) = - 2F_{\alpha}^3, \quad
\rho_{\alpha}(\P_3(h_2, h_2)) =  2E_{\alpha}^3,\quad
\rho_{\alpha}(\P_3(h_1, h_2)) =  H_{\alpha}^3,\\
&\rho_{\alpha}(e_1\otimes f_1\otimes h_1) = {\sqrt 2}iT_{\alpha}^1, \quad
\rho_{\alpha}(e_1\otimes f_1\otimes h_2) = {\sqrt 2}iT_{\alpha}^2, \\
&\rho_{\alpha}(e_1\otimes f_2\otimes h_1) = -{\sqrt 2}iT_{\alpha}^4, \quad
\rho_{\alpha}(e_1\otimes f_2\otimes h_2) = {\sqrt 2}iT_{\alpha}^3,\\
&\rho_{\alpha}(e_2\otimes f_1\otimes h_1) = {\sqrt 2}iD_{\alpha}^3, \quad
\rho_{\alpha}(e_2\otimes f_1\otimes h_2) = {\sqrt 2}iD_{\alpha}^4, \\
&\rho_{\alpha}(e_2\otimes f_2\otimes h_1) = -{\sqrt 2}iD_{\alpha}^2, \quad
\rho_{\alpha}(e_2\otimes f_2\otimes h_2) = {\sqrt 2}iD_{\alpha}^1. \\
\end{aligned}
\tag{19}
\end{equation*}
Thus $sp(\psi_i) \cong \hbox{Span} (E_{\alpha}^i, H_{\alpha}^i, F_{\alpha}^i)$
for $i = 1, 2, 3$.
$$\eqno\Q$$

\vfil\eject
\noindent
{\bf 3. Matrices over a Weyl algebra}

\vskip 0.5in

By definition, a Weyl algebra is
$$\W = \sum_{i\geq 0}\A d^i, \eqno (20)$$
 where $\A$ is an associative commutative algebra and
 $d:\A\rightarrow \A$ is a derivation of $\A$,
 with the relations
$$da = d(a) + ad, \quad a\in\A,\eqno (21)$$
see [13]. Set
$$\A = \C[t, t^{-1}], \hbox{ } d= {\partial\over{\partial t}}. \eqno (22)$$
Note that $\W$ is isomorphic to the  associative algebra $P_{\hbox{h} = 1}^+$.
Let $\hbox{End}(\W^{2|2})$ be the complex Lie superalgebra of $4\times 4$
matrices over  $\W$.

\vskip 0.1in
\noindent
{\bf Theorem 3.1.} For each
$\alpha \in \C$  there exists an embedding
$$\bar{\rho}_{\alpha}: \Gamma (2, -1 - \alpha, \alpha - 1)\rightarrow \hbox{End}(\W^{2|2})\eqno(23)$$
 given as follows:
\bigskip
\begin{equation*}
\begin{aligned}
&\bar{\rho}_{\alpha}(T_{\alpha}^{1}) = \left(\begin{array}{cc|cc}
0&0&t&0\\
0&0&0&0\\
\hline
0&0&0&0\\
0&t&0&0
\end{array}\right), \quad
\bar{\rho}_{\alpha}(T_{\alpha}^2) = \left(\begin{array}{cc|cc}
0&0&0&t\\
0&0&0&0\\
\hline
0&-t&0&0\\
0&0&0&0
\end{array}\right),\\
&\bar{\rho}_{\alpha}(T_{\alpha}^{3}) = \left(\begin{array}{cc|cc}
0&0&0&0\\
0&0&0&t\\
\hline
t&0&0&0\\
0&0&0&0
\end{array}\right), \quad
\bar{\rho}_{\alpha}(T_{\alpha}^4) = \left(\begin{array}{cc|cc}
0&0&0&0\\
0&0&-t&0\\
\hline
0&0&0&0\\
t&0&0&0
\end{array}\right),\\
&\bar{\rho}_{\alpha}(D_{\alpha}^{1}) = \left(\begin{array}{cc|cc}
0&0&0&0\\
0&0&0&d + \alpha t^{-1}\\
\hline
d&0&0&0\\
0&0&0&0
\end{array}\right), \quad
\bar{\rho}_{\alpha}(D_{\alpha}^2) = \left(\begin{array}{cc|cc}
0&0&0&0\\
0&0&-d - \alpha t^{-1}&0\\
\hline
0&0&0&0\\
d&0&0&0
\end{array}\right),\\
&\bar{\rho}_{\alpha}(D_{\alpha}^{3}) = \left(\begin{array}{cc|cc}
0&0&d+\alpha t^{-1}&0\\
0&0&0&0\\
\hline
0&0&0&0\\
0&d&0&0
\end{array}\right), \quad
\bar{\rho}_{\alpha}(D_{\alpha}^4) = \left(\begin{array}{cc|cc}
0&0&0&d+\alpha t^{-1}\\
0&0&0&0\\
\hline
0&-d&0&0\\
0&0&0&0
\end{array}\right),
\end{aligned}
\tag{24}
\end{equation*}
\begin{equation*}
\begin{aligned}
&\bar{\rho}_{\alpha}(E_{\alpha}^{1}) = t^2 1_{4|4}\\
&\bar{\rho}_{\alpha}(F_{\alpha}^1) = \left(\begin{array}{c|c}
 (d^2 + \alpha t^{-1}d)1_{2|2}&0\\
\hline
0&(d^2 + \alpha dt^{-1})1_{2|2}\\
\end{array}\right),\\
&\bar{\rho}_{\alpha}(H_{\alpha}^{1}) =
(td + {{1 +\alpha}\over 2})1_{4|4},
\end{aligned}
\end{equation*}
\begin{equation*}
\begin{aligned}
&\bar{\rho}_{\alpha}(E_{\alpha}^{2}) = \left(\begin{array}{cc|cc}
0&0&0&0\\
1&0&0&0\\
\hline
0&0&0&0\\
0&0&0&0
\end{array}\right),
\bar{\rho}_{\alpha}(F_{\alpha}^2) = \left(\begin{array}{cc|cc}
0&-1&0&0\\
0&0&0&0\\
\hline
0&0&0&0\\
0&0&0&0
\end{array}\right),
\bar{\rho}_{\alpha}(H_{\alpha}^{2}) = \left(\begin{array}{cc|cc}
-1&0&0&0\\
0&1&0&0\\
\hline
0&0&0&0\\
0&0&0&0
\end{array}\right),\\
&\bar{\rho}_{\alpha}(E_{\alpha}^{3}) = \left(\begin{array}{cc|cc}
0&0&0&0\\
0&0&0&0\\
\hline
0&0&0&1\\
0&0&0&0
\end{array}\right),
\bar{\rho}_{\alpha}(F_{\alpha}^3) = \left(\begin{array}{cc|cc}
0&0&0&0\\
0&0&0&0\\
\hline
0&0&0&0\\
0&0&1&0
\end{array}\right),
\bar{\rho}_{\alpha}(H_{\alpha}^{3}) = \left(\begin{array}{cc|cc}
0&0&0&0\\
0&0&0&0\\
\hline
0&0&1&0\\
0&0&0&-1
\end{array}\right).
\end{aligned}
\end{equation*}
{\it Proof.}
For each $\hbox{h}\in (0, 1]$ and each $\alpha\in\C$
there exists an  embedding
$$\rho_{\alpha, \hbox{h}}: \Gamma (2, -1 - \alpha, \alpha - 1)
\rightarrow P_{\hbox{h}}^+(4).\eqno (25)$$

\noindent
$\Gamma_{\alpha, \hbox{h}} = \rho_{\alpha, \hbox{h}}(\Gamma (2, -1 - \alpha, \alpha - 1))$
is spanned by the following elements:

\begin{equation*}
\begin{aligned}
&E_{\alpha,\hbox{h}}^1 =  t^2,\quad H_{\alpha,\hbox{h}}^1 = t\tau + {{\alpha + 1}\over 2}\hbox{h},\\
&F_{\alpha,\hbox{h}}^1 = \tau^2 - \alpha(2t^{-2}\xi_1\xi_2\eta_1\eta_2 +
t^{-2}(\xi_1\eta_1 + \xi_2\eta_2)\hbox{h} - t^{-1}\tau \hbox{h}),\\
&E_{\alpha,\hbox{h}}^2 = \xi_1\xi_2,\quad H_{\alpha,\hbox{h}}^2 =  \xi_1\eta_1 + \xi_2\eta_2 - \hbox{h},\quad
F_{\alpha,\hbox{h}}^2 = \eta_1\eta_2,\\
&E_{\alpha,\hbox{h}}^3 = \xi_1\eta_2,\quad H_{\alpha,\hbox{h}}^3 =  \xi_1\eta_1 - \xi_2\eta_2,\quad
F_{\alpha,\hbox{h}}^3 = \xi_2\eta_1,\\
&T_{\alpha,\hbox{h}}^1 =   t\eta_1,\quad
T_{\alpha,\hbox{h}}^2 =  t\eta_2,\quad
T_{\alpha,\hbox{h}}^3 =  t\xi_1,\quad
T_{\alpha,\hbox{h}}^4 =  t\xi_2,\\
&D_{\alpha,\hbox{h}}^1 =  \tau\xi_1 + \alpha t^{-1}\xi_1\xi_2\eta_2,\quad
D_{\alpha,\hbox{h}}^2 = \tau\xi_2 - \alpha t^{-1}\xi_1\xi_2\eta_1,\\
&D_{\alpha,\hbox{h}}^3 = \tau\eta_1 + \alpha t^{-1}\eta_1\eta_2\xi_2,\quad
D_{\alpha,\hbox{h}}^4 = \tau\eta_2 - \alpha t^{-1}\eta_1\eta_2\xi_1,\\
\end{aligned}
\tag{26}
\end{equation*}
and so
$$\hbox{lim}_{\hbox{h}\rightarrow 0}\Gamma_{\alpha, \hbox{h}} =
\Gamma_{\alpha}\subset P^+(4).\eqno (27)$$
We fix $\hbox{h} = 1$.
Let $V = \C[t, t^{-1}]\otimes\Lambda(\xi_1, \xi_2)$.
Define a representation of $\Gamma (2, -1 - \alpha, \alpha - 1)$ in $V$
according to the embedding ${\rho}_{\alpha,\hbox{h}=1}$. Namely,
$\xi_i$ is the operator of multiplication in $\Lambda(\xi_1, \xi_2)$,
$\eta_i$ is identified with $\partial_{\xi_i}$,
and $1\in P_{\hbox{h}=1}^+(4)$ acts by the identity operator.
Consider the following basis in $V$:

\begin{equation*}
\begin{aligned}
&v_m^0 = t^{m},\quad v_m^1 = t^{m}\xi_1,\\
&v_m^2 = t^{m}\xi_2, \quad v_m^{3} = t^{m}\xi_1\xi_2
\hbox{ for all }m\in\Z.
\end{aligned}
\tag{28}
\end{equation*}
Explicitly, the action of $\Gamma (2, -1 - \alpha, \alpha - 1)$ on $V$
is given as follows
\begin{equation*}
\begin{aligned}
&T_{\alpha}^1(v_m^3) = v_{m+1}^2,\quad
T_{\alpha}^1(v_m^{1}) = v_{m+1}^0,\quad
T_{\alpha}^2(v_m^3) = -v_{m+1}^1,\quad
T_{\alpha}^2(v_m^{2}) = v_{m+1}^0,\\
&T_{\alpha}^3(v_m^0) = v_{m+1}^1,\quad
T_{\alpha}^3(v_m^{2}) = v_{m+1}^{3},\quad
T_{\alpha}^4(v_m^0) = v_{m+1}^2,\quad
T_{\alpha}^4(v_m^{1}) = -v_{m+1}^{3},\\
&D_{\alpha}^1(v_m^0) = mv_{m-1}^1,\quad
D_{\alpha}^1(v_m^{2}) = (m + \alpha)v_{m-1}^{3},\\
&D_{\alpha}^2(v_m^0) = v_{m-1}^2,\quad
D_{\alpha}^2(v_m^{1}) = -(m + \alpha)v_{m-1}^{3},\\
&D_{\alpha}^3(v_m^{3}) = mv_{m-1}^2,\quad
D_{\alpha}^3(v_m^{1}) = (m + \alpha)v_{m-1}^{0},\\
&D_{\alpha}^4(v_m^{3}) = -mv_{m-1}^1,\quad
D_{\alpha}^4(v_m^{2}) = (m +  \alpha)v_{m-1}^{0},\\
&E_{\alpha}^1(v_m^0) = v_{m+2}^0,\quad
E_{\alpha}^1(v_m^3) = v_{m+2}^3,\quad
E_{\alpha}^1(v_m^1) = v_{m+2}^1,\quad
E_{\alpha}^1(v_m^2) = v_{m+2}^1,\\
&F_{\alpha}^1(v_m^0) = m(m - 1 + \alpha ) v_{m-2}^0,\quad
F_{\alpha}^1(v_m^3) = m(m - 1 + \alpha) v_{m-2}^3,\\
&F_{\alpha}^1(v_m^1) = (m + \alpha)(m - 1) v_{m-2}^1,\quad
F_{\alpha}^1(v_m^2) = (m + \alpha)(m - 1) v_{m-2}^2,\\
&H_{\alpha}^1(v_m^i) = (m + {{\alpha + 1}\over 2}) v_{m}^i,
\quad i = 0, 1, 2, 3,\\
&E_{\alpha}^2(v_m^0) =
v_{m}^{3},\quad
 F_{\alpha}^2(v_m^{3}) = -v_{m}^0,\quad
H_{\alpha}^2(v_m^0) = -v_{m}^0,\quad
H_{\alpha}^2(v_m^{3}) = v_{m}^{3},\\
&E_{\alpha}^3(v_m^2) = v_{m}^{1},\quad
 F_{\alpha}^3(v_m^{1}) = v_{m}^2,\quad
H_{\alpha}^3(v_m^1) = v_{m}^1,\quad H_{\alpha}^3(v_m^{2}) = -v_{m}^{2}.\\
\end{aligned}
\tag{29}
\end{equation*}
Thus we obtain the above-mentioned
embedding $\bar{\rho}_{\alpha}$ of $\Gamma (2, -1 - \alpha, \alpha - 1)$ into
$\hbox{End}(\W^{2|2})$.
$$\eqno\Q$$
{\bf Remark 3.2.} Note that the superalgebras $\Gamma (2, -1 - \alpha, \alpha - 1)$
are not simple, when $\alpha = 1$ or $-1$. Correspondingly,
we have the following realizations of $\p\s\l(2|2)$ as
a subsuperalgebra of $\hbox{End}(\W^{2|2})$.

\noindent
If $\alpha = 1$, then
\begin{equation*}
\begin{aligned}
&\hbox{Span} (E^i_{\alpha}, H^i_{\alpha}, F^i_{\alpha},  T^j_{\alpha}, D^j_{\alpha}
\hbox{ }|\hbox{ } i = 1, 2 \hbox{ and } j = 1, \ldots, 4) \cong \p\s\l(2|2),\\
&\Gamma (2, -2, 0)/\p\s\l(2|2) \cong \s\l(2).
\end{aligned}
\tag{30}
\end{equation*}
If $\alpha = - 1$, then
\begin{equation*}
\begin{aligned}
&\hbox{Span} (E^i_{\alpha}, H^i_{\alpha}, F^i_{\alpha},  T^j_{\alpha}, D^j_{\alpha}
\hbox{ } |\hbox{ } i = 1, 3 \hbox{ and } j = 1, \ldots, 4) \cong \p\s\l(2|2),\\
&\Gamma (2, 0, -2)/\p\s\l(2|2) \cong \s\l(2).
\end{aligned}
\tag{31}
\end{equation*}

\vskip 0.5in
\noindent

{\bf 4. Contractions of  $\Gamma (2, -1 - \alpha, \alpha - 1)$}

\vskip 0.5in

We denote the associative algebra $P_{\hbox{h}=1}$ by $\tilde{\W}$.
The product $A(t, \tau)\circ B(t, \tau)$ in $\tilde{\W}$ is defined as in (8) where $\hbox{h} = 1$.
Clearly, $\W\subset\tilde{\W}$.
In this section we consider a contraction $\Gamma$ of $\Gamma (2, -1 - \alpha, \alpha - 1)$ when $\alpha\rightarrow  1$
(or $\alpha\rightarrow  -1$), which is
a Lie superalgebra isomorphic to the universal central extension of $\p\s\l(2|2)$ (cf. [6]).
We obtain an embedding of $\Gamma$ into $4\times 4$ matrices over $\tilde{\W}$. In this realization  we essentially use pseudodifferential operators.
 This  allows us to construct
a three-parameter family of irreducible representations of the universal central extension of $\p\s\l(2|2)$
in $(2|2)$-dimensional complex superspace.

\vskip 0.1in
The nonzero commutation relations in $\Gamma (2, -1 - \alpha, \alpha - 1)$ are as follows:
\begin{equation*}
\begin{aligned}
&[T_{\alpha}^i, T_{\alpha}^j]= E_{\alpha}^1,\quad [D_{\alpha}^i, D_{\alpha}^j]= F_{\alpha}^1, \quad
[E_{\alpha}^1, D_{\alpha}^i]= -2T_{\alpha}^j,\\
&[E_{\alpha}^1, D_{\alpha}^j]= -2T_{\alpha}^i,\quad
[F_{\alpha}^1, T_{\alpha}^i]= 2D_{\alpha}^j, \quad [F_{\alpha}^1, T_{\alpha}^j]= 2D_{\alpha}^i,\\
&\quad\hbox{where } i=1, j=3 \hbox{ or } i=2, j=4;\\
&[E_{\alpha}^2, T_{\alpha}^i]= \mp T_{\alpha}^j, \quad [E_{\alpha}^2, D_{\alpha}^j]= \pm D_{\alpha}^i,\\
&[F_{\alpha}^2, T_{\alpha}^j]= \pm T_{\alpha}^i, \quad [F_{\alpha}^2, D_{\alpha}^i]= \mp D_{\alpha}^j,\\
&\quad\hbox{where } i=1, j=4 \hbox{ or } i=2, j=3;\\
&[E_{\alpha}^1, F_{\alpha}^1]= -4H_{\alpha}^1, \quad [E_{\alpha}^2, F_{\alpha}^2]= -H_{\alpha}^2, \\
&[H_{\alpha}^i, E_{\alpha}^i]= 2E_{\alpha}^i, \quad [H_{\alpha}^i, F_{\alpha}^i]= -2F_{\alpha}^i,
\quad i=1, 2,\\
&[H_{\alpha}^1, T_{\alpha}^i]= T_{\alpha}^i, \quad [H_{\alpha}^1, D_{\alpha}^i]= -D_{\alpha}^i,\quad i = 1, \ldots, 4\\
&[H_{\alpha}^2, T_{\alpha}^i]= \mp T_{\alpha}^i, \quad [H_{\alpha}^2, D_{\alpha}^i]= \pm D_{\alpha}^i,
\quad i = 1,2 \hbox{ or } i = 3, 4. \\
\end{aligned}
\tag{32}
\end{equation*}
\vskip 0.1in
\begin{equation*}
\begin{aligned}
&[E_{\alpha}^3, T_{\alpha}^i]= \mp T_{\alpha}^j, \quad [E_{\alpha}^3, D_{\alpha}^j]=  \pm D_{\alpha}^i,\\
&[F_{\alpha}^3, T_{\alpha}^j]= \mp T_{\alpha}^i, \quad
[F_{\alpha}^3, D_{\alpha}^i]=  \pm D_{\alpha}^j,\\
&\quad\hbox{where } i=1, j=2, \hbox{ or } i=4, j=3;\\
&[H_{\alpha}^3, T_{\alpha}^i]= \mp T_{\alpha}^i, \quad [H_{\alpha}^3, D_{\alpha}^i]= \pm D_{\alpha}^i,
\quad i = 1,4 \hbox{ or } i = 2, 3 \\
&[E_{\alpha}^3, F_{\alpha}^3]= -H_{\alpha}^3, \quad [H_{\alpha}^3, E_{\alpha}^3]= 2E_{\alpha}^3,
\quad [H_{\alpha}^3, F_{\alpha}^3]= -2F_{\alpha}^3.\\
\end{aligned}
\tag{33}
\end{equation*}
\vskip 0.1in

\begin{equation*}
\begin{aligned}
&[D_{\alpha}^i, T_{\alpha}^j]= \pm (1+\alpha)E_{\alpha}^2, \quad [T_{\alpha}^i, D_{\alpha}^j]= \mp (1+\alpha)F_{\alpha}^2,
\hbox{ where } i=1, j=4 \hbox{ or } i=2, j=3,\\
&[T_{\alpha}^i, D_{\alpha}^j]= \pm (\alpha - 1)E_{\alpha}^3, \quad [T_{\alpha}^j, D_{\alpha}^i]= \pm (\alpha - 1)F_{\alpha}^3,
\hbox{ where } i=3, j=4 \hbox{ or } i=2, j=1,\\
&[T_{\alpha}^1, D_{\alpha}^1]= H_{\alpha}^1 + ({{1+\alpha}\over 2})H_{\alpha}^2 + ({{1-\alpha}\over 2})H_{\alpha}^3,\quad
[T_{\alpha}^2, D_{\alpha}^2]= H_{\alpha}^1 + ({{1+\alpha}\over 2})H_{\alpha}^2 + ({{\alpha - 1}\over 2})H_{\alpha}^3\\
&[T_{\alpha}^3, D_{\alpha}^3]= H_{\alpha}^1 - ({{1+\alpha}\over 2})H_{\alpha}^2 + ({{\alpha - 1}\over 2})H_{\alpha}^3,\quad
[T_{\alpha}^4, D_{\alpha}^4]= H_{\alpha}^1 - ({{1+\alpha}\over 2})H_{\alpha}^2 + ({{1 - \alpha}\over 2})H_{\alpha}^3.\\
\end{aligned}
\tag{34}
\end{equation*}

\vskip 0.5in
\noindent
Let $$E_{\alpha}^3 = (\alpha - 1)^{-1}C_{+},\quad H_{\alpha}^3 = (\alpha - 1)^{-1}C,\quad F_{\alpha}^3 = (\alpha - 1)^{-1}C_{-}.
\eqno (35) $$
If $\alpha \rightarrow 1$, then $\Gamma (2, -1 - \alpha, \alpha - 1)$ contracts to the Lie superalgebra $\Gamma$,
so
$$\Gamma = \hbox{Span}(T_i, D_i, E_j, H_j, F_j, C, C_{+}, C_{-}) \quad\hbox{where } i = 1, \ldots, 4, j = 1, 2.\eqno (36)$$

\noindent
Note that the  nonzero commutation relations in $\Gamma$ are as in (32) and as follows:

\begin{equation*}
\begin{aligned}
&[T_{\alpha}^1, D_{\alpha}^4]= -2F_{\alpha}^2,\quad [T_{\alpha}^2, D_{\alpha}^3]= 2F_{\alpha}^2,\quad
[D_{\alpha}^1, T_{\alpha}^4]= 2E_{\alpha}^2, \quad [D_{\alpha}^2, T_{\alpha}^3]= -2E_{\alpha}^2, \\
&[T_{\alpha}^3, D_{\alpha}^4]= C_{+}, \quad [T_{\alpha}^4, D_{\alpha}^3]= C_{-},\quad
[T_{\alpha}^1, D_{\alpha}^2]= -C_{-}, \quad [T_{\alpha}^2, D_{\alpha}^1]= -C_{+},\\
&[T_{\alpha}^1, D_{\alpha}^1]= H_{\alpha}^1 + H_{\alpha}^2 - {C\over 2},\quad
[T_{\alpha}^2, D_{\alpha}^2]= H_{\alpha}^1 + H_{\alpha}^2 + {C\over 2},\\
&[T_{\alpha}^3, D_{\alpha}^3]= H_{\alpha}^1 - H_{\alpha}^2 + {C\over 2},\quad
[T_{\alpha}^4, D_{\alpha}^4]= H_{\alpha}^1 - H_{\alpha}^2 - {C\over 2}.\\
\end{aligned}
\tag{37}
\end{equation*}
{\bf Remark 4.1.} Let
$$E_{\alpha}^2 = (\alpha + 1)^{-1}C_{+},\quad H_{\alpha}^2 = (\alpha + 1)^{-1}C,\quad F_{\alpha}^2 = (\alpha + 1)^{-1}C_{-}.
\eqno (38)$$
If $\alpha \rightarrow -1$, then $\Gamma (2, -1 - \alpha, \alpha - 1)$ contracts to the Lie superalgebra
$$\hbox{Span}(T_i, D_i, E_j, H_j, F_j, C, C_{+}, C_{-}) \quad\hbox{where } i = 1, \ldots, 4, j = 1, 3,\eqno (39)$$
which is isomorphic to $\Gamma$.

\noindent
Let $a, b\in \C$ be such that $ab\not= 0,\pm 1$.

\noindent
{\bf Theorem 4.2.}
The Lie superalgebra $\Gamma$ is isomorphic to a  subsuperalgebra of $\hbox{End}(\tilde{\W}^{2|2})$
spanned by the following matrices:

\begin{equation*}
\begin{aligned}
&\tilde{C}_{+} = (t\tau){1_{2|2}}, \quad
\tilde{C} = {1_{2|2}}, \quad
\tilde{C}_{-} = ({\tau}^{-1}\circ t^{-1}){1_{2|2}},
\end{aligned}
\end{equation*}
\begin{equation*}
\begin{aligned}
&\tilde{E}_1 = \left(\begin{array}{cc|cc}
0&1&0&0\\
0&0&0&0\\
\hline
0&0&0&0\\
0&0&0&0
\end{array}\right),\quad
\tilde{F}_1 = \left(\begin{array}{cc|cc}
0&0&0&0\\
1&0&0&0\\
\hline
0&0&0&0\\
0&0&0&0
\end{array}\right),\quad
\tilde{H}_1 = \left(\begin{array}{cc|cc}
1&0&0&0\\
0&-1&0&0\\
\hline
0&0&0&0\\
0&0&0&0
\end{array}\right),
\end{aligned}
\end{equation*}
\begin{equation*}
\begin{aligned}
&\tilde{E}_2 = \left(\begin{array}{cc|cc}
0&0&0&0\\
0&0&0&0\\
\hline
0&0&0&0\\
0&0&1&0
\end{array}\right),\quad
\tilde{F}_2 = \left(\begin{array}{cc|cc}
0&0&0&0\\
0&0&0&0\\
\hline
0&0&0&1\\
0&0&0&0
\end{array}\right),\quad
\tilde{H}_2 = \left(\begin{array}{cc|cc}
0&0&0&0\\
0&0&0&0\\
\hline
0&0&1&0\\
0&0&0&-1
\end{array}\right),
\end{aligned}
\end{equation*}
\begin{equation*}
\begin{aligned}
&\tilde{T}_3 = \left(\begin{array}{cc|cc}
0&0&t\tau&0\\
0&0&0&0\\
\hline
0&0&0&0\\
0&a&0&0
\end{array}\right),\quad
\tilde{T}_2 = \left(\begin{array}{cc|cc}
0&0&0&-t\tau\\
0&0&0&0\\
\hline
0&a&0&0\\
0&0&0&0
\end{array}\right),
\end{aligned}
\end{equation*}
\begin{equation*}
\begin{aligned}
&\tilde{D}_4 = \left(\begin{array}{cc|cc}
0&0&0&0\\
0&0&0&t\tau\\
\hline
a&0&0&0\\
0&0&0&0
\end{array}\right),\quad
\tilde{D}_1 = \left(\begin{array}{cc|cc}
0&0&0&0\\
0&0&-t\tau&0\\
\hline
0&0&0&0\\
a&0&0&0
\end{array}\right).
\end{aligned}
\end{equation*}
\begin{equation*}
\begin{aligned}
&\tilde{T}_1 = \left(\begin{array}{cc|cc}
0&0&0&b\\
0&0&0&0\\
\hline
0&\tau^{-1}\circ t^{-1}&0&0\\
0&0&0&0
\end{array}\right),\quad
\tilde{T}_4 = \left(\begin{array}{cc|cc}
0&0&b&0\\
0&0&0&0\\
\hline
0&0&0&0\\
0&-\tau^{-1}\circ t^{-1}&0&0
\end{array}\right),
\end{aligned}
\end{equation*}
\begin{equation*}
\begin{aligned}
&\tilde{D}_2 = \left(\begin{array}{cc|cc}
0&0&0&0\\
0&0&b&0\\
\hline
0&0&0&0\\
\tau^{-1}\circ t^{-1}&0&0&0
\end{array}\right),\quad
\tilde{D}_3 = \left(\begin{array}{cc|cc}
0&0&0&0\\
0&0&0&b\\
\hline
-\tau^{-1}\circ t^{-1}&0&0&0\\
0&0&0&0
\end{array}\right).
\end{aligned}
\tag{40}
\end{equation*}
{\it Proof.} Let $q = \sqrt{{2\over {1 + ab}}}$

An embedding
$$\varphi: \Gamma\longrightarrow \hbox{End}(\tilde{\W}^{2|2})\eqno (41)$$
is given as follows:
\begin{equation*}
\begin{aligned}
&\varphi(T_1) = -qi\tilde{T_1}, \quad \varphi(T_2) = q\tilde{T_2}, \quad
\varphi(T_3) = qi\tilde{T_3}, \quad \varphi(T_4) = q\tilde{T_4}, \\
&\varphi(D_1) = qi\tilde{D_1}, \quad \varphi(D_2) = -q\tilde{D_2}, \quad
\varphi(D_3) = qi\tilde{D_3}, \quad \varphi(D_4) = q\tilde{D_4}, \\
&\varphi(E_1) = 2\tilde{E_1}, \quad \varphi(H_1) = \tilde{H_1}, \quad \varphi(F_1) = -2\tilde{F_1},\\
&\varphi(E_2) = i\tilde{E_2}, \quad \varphi(H_2) = -\tilde{H_2}, \quad \varphi(F_2) = i\tilde{F_2},\\
&\varphi(C) = q^2(1 - ab)\tilde{C}, \quad \varphi(C_{+}) = q^2ai\tilde{C_{+}}, \quad \varphi(C_{-}) = -q^2bi\tilde{C}_{-}\\
\end{aligned}
\tag{42}
\end{equation*}
$$\eqno\Q$$
Recall that
$$\s\l(2|2) = \lbrace
\left(\begin{array}{c|c}
A&B\\
\hline
C&D\\
\end{array}\right) \hbox{ }|\hbox{ }A, B, C, D \hbox{ are } 2\times 2\hbox{ matrices over }\C  \hbox{ such that } tr A = tr D\rbrace, \eqno (43)$$
$$\p\s\l(2|2) = \s\l(2|2)/<1_{2|2}>,\eqno (44)$$
see [1].
Let ${B}_{ij}$, ${C}_{ij}$ be elementary $2\times 2$ matrices.
Let $\tilde{B}_{ij}, \tilde{C}_{ij}$ and  $\tilde{E}_i, \tilde{F}_i, \tilde{H}_i$
(see (40)), where $i = 1, 2$, be the corresponding  elements of $\p\s\l(2|2)$.

\noindent
Recall that a central extension of a Lie superalgebra $\g$ is a pair $(U, \psi)$ consisting of a Lie superalgebra $U$ and
a surjective homomorphism $\psi: U\rightarrow \g$ such that $[\hbox{Ker}\psi, U] = 0$.

\noindent
A central extension $(U, \psi)$  of a Lie superalgebra $\g$ is universal if the following conditions hold (see [2]):
\begin{equation*}
\begin{aligned}
&(1) \quad [\g, \g] = \g,\\
&(2)\quad \hbox{for any central extension }
(W, \phi)\hbox{ of } \g \hbox{ there exists a homomorphism}\\
& \qquad\nu: U\rightarrow W \hbox{ such that } \phi\circ\nu = \psi.\\
\end{aligned}
\end{equation*}

\noindent
{\bf Remark 4.3.} Let $\g = \p\s\l(2|2)$. Note that $\hbox{dim}H^2(\g, \C) = 3$ [2].

\noindent
Let $\c= \hbox{Span}(c_+, c, c_-)\cong \C^3$.
The universal central extension of $\g$ is given by an exact sequence of Lie superalgebras
$$0\longrightarrow\c
\longrightarrow\hat{\g}\longrightarrow\g
\longrightarrow 0, \eqno (45)$$
so that
$$[(g_1, x_1), (g_2, x_2)] = ([g_1, g_2], \hbox{f}(g_1, g_2))\eqno (46)$$
where $g_1, g_2 \in \g$  and $x_1, x_2\in \hbox{Span}(c_+, c, c_-)$.
The corresponding $2$--cocycle $\hbox{f}$ is given as follows:
\begin{equation*}
\begin{aligned}
&\hbox{f}(\tilde{B}_{12}, \tilde{C}_{21}) = \hbox{f}(\tilde{C}_{12}, \tilde{B}_{21}) =
\hbox{f}(\tilde{C}_{22}, \tilde{B}_{22}) = \hbox{f}(\tilde{B}_{11}, \tilde{C}_{11}) = c,\\
&\hbox{f}(\tilde{C}_{22}, \tilde{C}_{11}) = -\hbox{f}(\tilde{C}_{12}, \tilde{C}_{21}) = c_+,\\
&\hbox{f}(\tilde{B}_{12}, \tilde{B}_{21}) = - \hbox{f}(\tilde{B}_{11}, \tilde{B}_{22}) = c_-.\\
\end{aligned}
\tag{47}
\end{equation*}
Let $\omega\in C^2(\g, \C)$, and let
$X_1, X_2, X_3\in {\g}$ be such that $p(X_1) = \bar{0}$, $p(X_2) = p(X_3) = \bar{1}$.
Then
$$d\omega(X_1, X_2, X_3) = -\omega([X_1, X_2], X_3) -\omega([X_1, X_3], X_2) -\omega([X_2, X_3], X_1).\eqno (48)$$
One can check that $d\hbox{f} = 0$; hence $\hbox{f}\in Z^2(\g, \C)$. On the other hand, if $\omega\in C^1(\g, \C)$ and
$p(X_1) = p(X_2) = \bar{1}$, then $d\omega(X_1, X_2) = -\omega([x_1, x_2])$. Since $[\tilde{C}_{22}, \tilde{C}_{11}] = 0$,
then $\hbox{f}\not\in B^2(\g, \C)$. Hence $\hbox{f}\in H^2(\g, \C)$, and it is easy to see that the conditions (1) and (2)
hold.

\noindent
{\bf Theorem 4.4.}
There exists a three-parameter family of irreducible representations  of the universal central extension
$\hat{\p\s\l}(2|2)$
in $(2|2)$--dimensional superspace depending on $\lambda, a, b \in \C$ such that $\lambda\not= 0$ and $ab \not= -1$.

\noindent
{\it Proof.}
Let $\lambda \not=0$.  Let  $V^{\lambda}$ be a complex $(2|2)$--dimensional superspace spanned by even vectors
$(t^{\lambda}, 0, 0, 0), (0, t^{\lambda}, 0, 0)$ and odd vectors $(0,0, t^{\lambda}, 0), (0, 0, 0, t^{\lambda})$. The natural action of $\varphi(\Gamma)$ on  $V^{\lambda}$
is given by the matrices $\tilde{E}_i, \tilde{F}_i, \tilde{H}_i$, where $i = 1, 2$, the matrices
$\tilde{C}_+^{\lambda} = (\lambda) 1_{2|2}$, $\tilde{C}_-^{\lambda} = ({\lambda}^{-1}) 1_{2|2}$ and
$\tilde{C} = 1_{2|2}$, and the matrices
$\tilde{T}_i^{\lambda}$, $\tilde{D}_i^{\lambda}$, where $i = 1, \ldots 4$:

\begin{equation*}
\begin{aligned}
&{\tilde{T}}_3^{\lambda} = \left(\begin{array}{cc|cc}
0&0&\lambda&0\\
0&0&0&0\\
\hline
0&0&0&0\\
0&a&0&0
\end{array}\right),\quad
{\tilde{T}}_2^{\lambda} = \left(\begin{array}{cc|cc}
0&0&0&-\lambda\\
0&0&0&0\\
\hline
0&a&0&0\\
0&0&0&0
\end{array}\right),
\end{aligned}
\end{equation*}
\begin{equation*}
\begin{aligned}
&{\tilde{D}}_4^{\lambda} = \left(\begin{array}{cc|cc}
0&0&0&0\\
0&0&0&\lambda\\
\hline
a&0&0&0\\
0&0&0&0
\end{array}\right),\quad
{\tilde{D}}_1^{\lambda} = \left(\begin{array}{cc|cc}
0&0&0&0\\
0&0&-\lambda&0\\
\hline
0&0&0&0\\
a&0&0&0
\end{array}\right).
\end{aligned}
\end{equation*}
\begin{equation*}
\begin{aligned}
&{\tilde{T}}_1^{\lambda} = \left(\begin{array}{cc|cc}
0&0&0&b\\
0&0&0&0\\
\hline
0&{\lambda}^{-1}&0&0\\
0&0&0&0
\end{array}\right),\quad
{\tilde{T}}_4^{\lambda} = \left(\begin{array}{cc|cc}
0&0&b&0\\
0&0&0&0\\
\hline
0&0&0&0\\
0&-{\lambda}^{-1}&0&0
\end{array}\right),
\end{aligned}
\end{equation*}
\begin{equation*}
\begin{aligned}
&{\tilde{D}}_2^{\lambda} = \left(\begin{array}{cc|cc}
0&0&0&0\\
0&0&b&0\\
\hline
0&0&0&0\\
{\lambda}^{-1}&0&0&0
\end{array}\right),\quad
{\tilde{D}}_3^{\lambda} = \left(\begin{array}{cc|cc}
0&0&0&0\\
0&0&0&b\\
\hline
-{\lambda}^{-1}&0&0&0\\
0&0&0&0
\end{array}\right).
\end{aligned}
\tag{49}
\end{equation*}
These matrices span a superalgebra isomorphic to a central extension of ${\p\s\l}(2|2)$.
Assume that $ab \not= -1$.
Let $s = {1\over {\sqrt{ab+1}}}$.
Explicitly a family of irreducible representations $\theta$ of the universal central extension $\hat{\p\s\l}(2|2)$ in $(2|2)$--dimensional complex superspace
is given as follows:
\begin{equation*}
\begin{aligned}
&\theta (\tilde{C}_{11}) = s{\tilde{D}}_4^{\lambda}, \quad
\theta (\tilde{C}_{12}) = s{\tilde{T}}_2^{\lambda}, \quad
\theta (\tilde{C}_{22}) = s{\tilde{T}}_3^{\lambda}, \quad
\theta (\tilde{C}_{21}) = s{\tilde{D}}_1^{\lambda}, \\
&\theta (\tilde{B}_{11}) = s{\tilde{T}}_4^{\lambda}, \quad
\theta (\tilde{B}_{12}) = s{\tilde{T}}_1^{\lambda}, \quad
\theta (\tilde{B}_{22}) = s{\tilde{D}}_3^{\lambda}, \quad
\theta (\tilde{B}_{21}) = s{\tilde{D}}_2^{\lambda}, \\
&\theta (\tilde{E}_{i}) =  \tilde{E}_{i},\quad
\theta (\tilde{F}_{i}) = \tilde{F}_{i}, \quad
\theta (\tilde{H}_{i}) = \tilde{H}_{i}\quad i = 1, 2,\\
&\theta ({c}) =  {{s^2(ab - 1)}\over 2} 1_{2|2},\quad
\theta ({c}_{+}) = (s^2\lambda a) 1_{2|2}, \quad
\theta ({c}_{-}) = (s^2{\lambda}^{-1}b)1_{2|2}.
\end{aligned}
\tag{50}
\end{equation*}
$$\eqno\Q$$

\vskip 0.3in
\noindent
{\bf Acknowledgements}
\vskip 0.3in

The author would like to thank the Max-Planck-Institut f\"ur Mathematik in Bonn
for hospitality and support in the summer of 2008.

The author is grateful to the referee for very valuable remarks and for references
on Jordan superalgebras.
\vskip 0.5in
\noindent
{\bf References}

\vskip 0.3in

\begin {itemize}

\font\red=cmbsy10
\def\~{\hbox{\red\char'0016}}

\item[{[1]}]
V. G. Kac,
Lie superalgebras,
Adv. Math. 26 (1977) 8--96.

\item[{[2]}] K. Iohara , Y. Koga,
Central extensions of Lie superalgebras,
Comment. Math. Helv. 76 (2001) 110--154.

\item[{[3]}] N. Beisert,
The $\s\u (2|2)$ dynamic S-matrix,
Adv. Theor. Math. Phys. 12 (2008) 945--979.
e-print arXiv:hep-th/0511082.

\item[{[4]}] N. Beisert, F. Spill,
The classical $r$-matrix of AdS/CFT and its Lie bialgebra structure,
Comm. Math. Phys. 285 (2009) 537--565.
e-print arXiv:0708.1762.

\item[{[5]}] I. Heckenberger, F. Spill, A. Torrelli, and H. Yamane,
 Drinfeld second realization of the affine superalgebras of $D^{(1)}(2, 1; x)$ via the Weyl groupoid, in:
Combinatorial Representation Theory and Related Topics, RIMS K{\^o}ky{\^u}roku Bessatsu, B8, Res. Inst. Math. Sci. (RIMS), Kyoto, 2008, pp. 171--216.
 e-print arXiv:0705.1071.

\item[{[6]}] T. Matsumoto and S. Moriyama,
An exceptional algebraic origin of the AdS/CFT Yangian symmetry,
J. High Energy Phys. 022 (4) (2008), 19 pp.
e-print arXiv:0803.1212.

\item[{[7]}] M. Scheunert, The Theory of Lie Superalgebras, in:
Lecture Notes in Mathematics vol. 716,  Springer, Berlin, 1979.

\item[{[8]}]
E. Poletaeva,
Embedding of the Lie superalgebra $D(2, 1 ; \alpha)$
into the Lie superalgebra of pseudodifferential symbols on $S^{1|2}$,
J. Math. Phys. 48 (2007) 103504, 17 pp.
e-print arXiv:0709.0083.

\item[{[9]}]
E. Poletaeva,
Deforming the Lie superalgebra $D(2, 1 ; \alpha)$
inside the superconformal algebra $K'(4)$.
Journal of Mathematical Sciences  161, no.1 (2009) 130--142.

\item[{[10]}] E. Poletaeva,
On matrix realizations of the contact superconformal algebra
$\hat{K}'(4)$ and the exceptional $N = 6$ superconformal algebra,
DCDIS A Supplement, Adv.  Dynam. Syst. 14 (S2) (2007) 285--289. e-print
arXiv:0707.3097.

\item[{[11]}] E. Poletaeva,
Matrix realizations of exceptional superconformal algebras.
J. Geom.  Phys. 58 (2008) 761--772.
e-print arXiv:0801.4398.

\item[{[12]}]
V. G. Kac,
Classification of supersymmetries, in:
Proceedings of the International Congress of Mathematicians, Beijing, 2002 (Higher Education Press, Beijing, 2002), Vol. I,  319--344.

\item[{[13]}] C. Martinez, E. I. Zelmanov,
Lie superalgebras graded by $P(n)$ and $Q(n)$,
Proc. Natl. Acad. Sci. USA 100 (2003)  8130--8137.

 \item[{[14]}]
C. Martinez, I. Shestakov and E. Zelmanov,
Jordan superalgebras defined by brackets,
J. London Math. Soc. (2) 64 (2001) 357---368.

\item[{[15]}]
I. L. Kantor,
Jordan and Lie superalgebras defined by Poisson algebras, in:
Algebra and Analysis. Tomsk, 1989. Amer. Math. Soc. Transl. Ser. 2 151 (1992) 55--79.

\item[{[16]}]
D. King and K. McCrimmon,
The Kantor  construction of Jordan superalgebras,
Comm. Algebra 20 (1992) 109--126.

\item[{[17]}] S. J. Cheng and V. G. Kac,
A new $N = 6$ superconformal algebra,
Comm. Math. Phys. 186 (1997) 219--231.

\item[{[18]}]
K. McCrimmon,
Speciality and non-speciality of two Jordan superalgebras,
J. Algebra 149 (1992) 326--351.

\item[{[19]}]
V. Ovsienko and C. Roger,
Deforming the Lie algebra of vector fields on $S^1$
inside the Poisson algebra on $\dot{T}^*S^1$,
Comm. Math. Phys. 198 (1998) 97--110.

\item[{[20]}]
V. Ovsienko and C. Roger,
Deforming the Lie algebra of vector fields on $S^1$
inside the Lie algebra of pseudodifferential symbols on $S^1$,
Amer. Math. Soc. Transl.  194 (1999) 211--226.

\end{itemize}

\end{document}